# POISSON–VORONOI APPROXIMATION

By Matthias Heveling[1] and Matthias Reitzner

*Vienna University of Technology*

Let $X$ be a Poisson point process and $K \subset \mathbb{R}^d$ a measurable set. Construct the Voronoi cells of all points $x \in X$ with respect to $X$, and denote by $v_X(K)$ the union of all Voronoi cells with nucleus in $K$. For $K$ a compact convex set the expectation of the volume difference $V(v_X(K)) - V(K)$ and the symmetric difference $V(v_X(K) \triangle K)$ is computed. Precise estimates for the variance of both quantities are obtained which follow from a new jackknife inequality for the variance of functionals of a Poisson point process. Concentration inequalities for both quantities are proved using Azuma's inequality.

**1. Introduction.** Let $X$ be a stationary Poisson point process of intensity $\lambda$. Denote by $v_X(x)$ the (random) Voronoi cell of $x$ with respect to the point set $X \cup \{x\}$, that is,

$$v_X(x) = \{z \in \mathbb{R}^d : \|z - x\| \leq \|z - y\| \text{ for all } y \in X\}.$$

We call $x$ the nucleus of the Voronoi cell $v_X(x)$. The set of all Voronoi cells $v_X(x)$, $x \in X$, is the Poisson–Voronoi tessellation of $\mathbb{R}^d$. For a given set $A \subset \mathbb{R}^d$ we consider the *Poisson–Voronoi approximation* $v_X(A)$ of $A$ which consists of all Voronoi cells with nucleus in $A$,

$$v_X(A) := \bigcup_{x \in X \cap A} v_X(x).$$

The set $v_X(A)$ is a random approximation of $A$. In this paper, we discuss the quality of this approximation for a convex set $A$. In particular, we are interested in the convergence of $v_X(A)$ to $A$ when the intensity of the Poisson process tends to infinity.

Received July 2007; revised July 2008.
[1]Supported by the European Network PHD, MCRN-511953.

The general problem whether $v_X(A)$ approximates $A$ for "complicated" sets $A$ seems to be difficult; only partial answers are available (see Khmaladze and Toronjadze [7] and Penrose [14]).

Here we concentrate on the case of a compact convex set $K$ with nonempty interior, and its approximation $v_X(K)$, where we derive precise estimates for the volume $V(v_X(K))$ and the volume of the symmetric difference of $K$ and $v_X(K)$.

THEOREM 1. *Let $X$ be a stationary Poisson point process of intensity $\lambda$. If $K$ is a convex set of volume $V(K)$ and surface area $S(K)$, then*

$$\mathbb{E} V(v_X(K)) = V(K)$$

*and*

$$\operatorname{Var} V(v_X(K)) \leq c_d \lambda^{-1-1/d} S(K)$$

*with a constant $c_d$ independent of $K$ and $\lambda$. Further there are constants $c(K), \Lambda(K)$ such that*

$$\mathbb{P}(|V(v_X(K)) - V(K)| \geq t\sqrt{\lambda^{-1-1/d} S(K)})$$
$$\leq c(K) e^{-c'_d t^2 (k \ln \lambda)^{-1-1/d}} + 16\sqrt{d} \lambda^{-k+1} S(K)$$

*with $c'_d = 2^{-4} 3^{-2d} d^{-d-1/2}$ for $\lambda \geq \Lambda(K)$ and $k \geq 2$.*

Our proof of the theorem relies on the stationarity of the process. Its first part can be generalized to nonstationary Poisson processes with an absolute continuous intensity measure with density function $\lambda f(x) > 0$ (with respect to Lebesgue measure). In that case the volume on the right-hand side of the formula is replaced by the integral of the density function $f$ of the intensity measure over $K$. Observe that if this density vanishes on a set of positive measure, then we do not even have $\mathbb{E} V(v_X(K)) \to V(K)$ for certain sets $K$. The second part of the theorem concerning the variance and the tail probability is also subject to generalization and holds for nonstationary Poisson processes with bounded density function. The present form of the theorem, however, gives fastest access to our asymptotic results.

Since the expectation of the volume of the Poisson–Voronoi approximation coincides with the volume of the convex set itself, it is natural to ask for the volume of the symmetric difference $K \triangle v_X(K) = (K \setminus v_X(K)) \cup (v_X(K) \setminus K)$. It is known that for any bounded Borel set $A \subset \mathbb{R}^d$, one has

$$V(A \triangle v_X(A)) \to 0$$

almost surely as $\lambda \to \infty$. This was proved for $d = 1$ in Khmaladze and Toronjadze [7], and by Einmahl and Khmaladze [3] for any bounded Borel set $A \subset \mathbb{R}^d$ with $V(A_\varepsilon) \to V(A)$ for $\varepsilon \to 0$, where $A_\varepsilon = A + \varepsilon B^d$ is the Minkowski sum of $A$ and the ball $\varepsilon B^d$. The general case was proved by Penrose [14].

Here we concentrate on rates of convergence for convex sets $A$. We derive an asymptotically precise formula for the expectation and estimates for the variance and the tails. We believe that the estimates for the variance in Theorems 1 and 2 are best possible up to the choice of the constant. Denote by $\kappa_d$ the volume of the unit ball in $\mathbb{R}^d$, and by $\Gamma(\cdot)$ the Gamma function.

THEOREM 2. *Let $X$ be a stationary Poisson point process of intensity $\lambda$. If $K$ is a convex set of volume $V(K)$ and surface area $S(K)$, then*

$$(1) \qquad \mathbb{E}V(K \triangle v_X(K)) = c_\mathbb{E} \lambda^{-1/d} S(K)(1 + O(\lambda^{-1/d}))$$

*with $c_\mathbb{E} = \frac{2}{d}\kappa_d^{-1/d}\kappa_{d-1}\Gamma(\frac{1}{d})$. And*

$$\operatorname{Var} V(K \triangle v_X(K)) \leq c_d \lambda^{-1-1/d} S(K)$$

*with a constant $c_d$ independent of $K$ and $\lambda$. Further there are constants $c(K), \Lambda(K)$ such that*

$$\mathbb{P}(|V(K \triangle v_X(K)) - \mathbb{E}V(K \triangle v_X(K))| \geq t\sqrt{\lambda^{-1-1/d}S(K)})$$
$$\leq c(K)e^{-c'_d t^2 (k \ln \lambda)^{-1-1/d}} + 16\sqrt{d}\lambda^{-k+1}S(K)$$

*with $c'_d = 2^{-4}3^{-2d}d^{-d-1/2}$ for $\lambda \geq \Lambda(K)$ and $k \geq 2$.*

REMARK 1. A precise estimate for the error term in (1) is given in Section 5, Theorem 6.

REMARK 2. In both theorems the estimates for the tail probabilities are stated for $\lambda \geq \Lambda(K)$ sufficiently large. Set $r = 4\sqrt{d}(k\lambda^{-1}\ln\lambda)^{1/d}$. Then $\lambda$ is sufficiently large, if $r \leq 4$ and the volume of the parallel set $V(\partial K + rB^d)$ is bounded by $4rS(K)$.

REMARK 3. Since in both cases, the results for the expectation and variance only depend on the volume and surface area of $K$, they possibly hold for more general classes of sets. Yet our methods of proof make essential use of the convexity of $K$ (in particular Hadwiger's characterization theorem; see Section 2). In view of applications it would be of interest to extend our results to sets in the convex ring.

REMARK 4. Jeulin posed the following interesting problem: Is it better to approximate $V(K)$ by a single realization $v_X(K)$ where $X$ is a Poisson point process of intensity $\lambda = k\lambda_0$, or by the mean value of $k$ realizations of $v_{X_i}(K)$ where $X_i$, $i = 1, \ldots, k$, are independent Poisson point processes of intensity $\lambda_0$? Both estimates are unbiased. Yet the approximation by a single estimate should be preferred, since $\mathrm{Var}(v_X(K))$ is of order $k^{-1-1/d}\lambda_0^{-1-1/d}$ whereas $\mathrm{Var}(\frac{1}{k}\sum_{i=1}^{k} v_{X_i}(K))$ is of order $k^{-1}\lambda_0^{-1-1/d}$.

Our results have applications in nonparametric statistics (see Einmahl and Khmaladze [3], Section 3) and image analysis (reconstructing an image from its intersection with a Poisson point process; see [7]). Also the connection to quantization problems is obvious; see, for example, Chapter 9 in the book of Graf and Luschgy [4] which gives an excellent introduction into this topic. Quantization problems are connected to problems of numerical integration: approximate the volume $V_d(A) = \int_A dx$ of a set $A$ using its intersection with a point process $X \cap A$. As shown in Theorem 1 the volume of the Poisson–Voronoi approximation $v_X(A)$ is an unbiased estimator for $V_d(A)$, even for Borel sets $A$ which is pointed out in Section 4. An estimate for the volume of $A$ is also obtained if the number of points $X(A)$ that fall into $A$ is counted instead of the volume of the Poisson–Voronoi approximation. By the definition of $X$ one has

$$\mathbb{E} X(A) = \mathrm{Var}\, X(A) = \lambda V(A).$$

When the variation coefficient is regarded as a measure for the quality of the respective volume estimators, then using Poisson–Voronoi approximation is more precise than counting points of the Poisson point process—at least for convex sets $A$ and stationary Poisson point processes. It would be of interest to show that this is a general principle even for arbitrary point processes, for example for random lattices $X$.

An interesting open problem is to measure the quality of approximation of $K$ by $v_X(K)$ using the Hausdorff distance between both sets. We are not aware of any results in this direction.

For basic facts from integral geometry, stochastic geometry and Voronoi tessellations which are not explained in the following, we refer the reader to [13, 16, 17, 18]. The employed notions and results from the theory of convex bodies are found in [15].

We work in $d$-dimensional Euclidean vector space $\mathbb{R}^d$, with norm $d(x,y) = \|x-y\|$, and for closed sets $K \subset \mathbb{R}^d$ distance $d(x,K) = \min(d(x,y);\ y \in K)$. Its unit ball, $\{x \in \mathbb{R}^d : |x\| \leq 1\}$, is denoted by $B^d$, and $S^{d-1}$ is the unit sphere. The space of convex bodies (nonempty, compact, convex subsets) in $\mathbb{R}^d$ is denoted by $\mathcal{K}^d$ and the space of locally finite point sets in $\mathbb{R}^d$ is denoted by $\mathbf{N}$.

For a stationary Poisson point process, as usual, $X$ denotes the simple counting measure as well as its support, that is, $X(A)$ and $\mathrm{card}(X \cap A)$ have the same meaning. Its intensity measure $\Theta = \mathbb{E} X(\cdot)$ ($\mathbb{E}$ denotes mathematical expectation) is given by

$$\mathbb{E} X(\cdot) = \lambda \int_{\mathbb{R}^d} \mathbf{1}\{x \in \cdot\}\, dx = \lambda V(\cdot).$$

**2. Valuations and Delaunay triangulations.** A major tool for proving our results is valuations. A functional $\mu:\mathcal{K}^d\to\mathbb{R}$ is called a *valuation*, if for every $K,L\in\mathcal{K}^d$ with $K\cup L\in\mathcal{K}^d$,

$$\mu(K\cup L)+\mu(K\cap L)=\mu(K)+\mu(L)$$

holds. A monotone valuation satisfies $\mu(K)\leq\mu(L)$ if $K\subset L$. Valuations play an important role in convex geometry and integral geometry; for further references see [9, 11] and [12]. One of the most important results in this field is the following characterization theorem by Hadwiger [5]:

*A functional $\mu:\mathcal{K}^d\to\mathbb{R}$ is a monotone and rigid motion-invariant valuation if and only if there are constants $c_0, c_1,\ldots,c_d$ (uniquely determined by $\mu$) such that*

$$\mu(K)=c_dV_d(K)+\cdots+c_0V_0(K)$$

*for every $K\in\mathcal{K}^d$.*

Here $V_0(K),\ldots,V_d(K)$ are the quermassintegrals of $K$. In particular, $V_d(K)$ equals the volume $V(K)$, $2V_{d-1}(K)$ is the surface area $S(K)$, and $V_0(K)$ is a multiple of the Euler characteristic. For a modern proof of this theorem, see Klain [8].

In the following sections a valuation depending on the Delaunay mosaic of $X$ turns out to be of importance. To this end denote by $\mathcal{E}_X$ the edges of the Delaunay mosaic of $X$, that is, those segments $[x,y]$ with $x,y\in X$ such that $v_X(x)\cap v_X(y)$ is a facet of $v_X(x)$ and $v_X(y)$. Set

$$n_K[x,y]=2\mathbf{1}([x,y]\cap K\neq\varnothing)-\mathbf{1}(x\in K)-\mathbf{1}(y\in K),$$

where $\mathbf{1}(\cdot)$ denotes the indicator function. Observe that for points $x,y\in X$, with probability 1, $n_K[x,y]$ is the number of connected components of $[x,y]\setminus K$ if $[x,y]$ meets $K$, and $n_K[x,y]=0$ otherwise.

THEOREM 3. *Let $f:\mathbf{N}\times\mathbb{R}^d\to[0,\infty)$ be measurable. Define a functional $\mu:\mathcal{K}^d\to\mathbb{R}$ by $\mu(K):=\mathbb{E}\sum_{[x,y]\in\mathcal{E}_X}(f(X,x)+f(X,y))n_K[x,y]$. Then*

(2) $$\mu(K)=c_f(\lambda)S(K),$$

*where $S(K)$ is the surface area of $K$. Moreover, if $f$ has the scaling property $f(t\varphi,tx)=t^\alpha f(\varphi,x)$, then there is a constant $c_f$ such that*

$$c_f(\lambda)=c_f\lambda^{(d-\alpha-1)/d}.$$

PROOF. First we will show that $\mu$ can be expressed as the difference of two auxiliary monotone valuations. We define

$$\nu_1(K):=\mathbb{E}\sum_{[x,y]\in\mathcal{E}_X}\mathbf{1}([x,y]\cap K\neq\varnothing)(f(X,x)+f(X,y)).$$

Elementary considerations yield that for any two convex bodies $K$ and $L$ such that $K \cup L$ is convex we have

$$\mathbf{1}([x,y] \cap K \neq \varnothing) + \mathbf{1}([x,y] \cap L \neq \varnothing)$$
$$= \mathbf{1}([x,y] \cap (K \cup L) \neq \varnothing) + \mathbf{1}([x,y] \cap (K \cap L) \neq \varnothing).$$

Hence $\nu_1$ is a valuation and it follows directly from the definition that $\nu_1$ is also monotone and rigid motion-invariant. We define a second functional $\nu_2$ by

$$\nu_2(K) := \mathbb{E} \sum_{[x,y] \in \mathcal{E}_X} (\mathbf{1}(x \in K) + \mathbf{1}(y \in K))(f(X,x) + f(X,y)).$$

As above for $\nu_1$ it is easily seen that $\nu_2$ is also a monotone and rigid motion-invariant valuation. Then Hadwiger's theorem yields that both $\nu_1$ and $\nu_2$ can be written as a linear combination of the Minkowski functionals. Moreover, since $\mu = 2\nu_1 - \nu_2$, we can write $\mu(K) = \sum_{i=0}^{d} c_i(\lambda) V_i(K)$. If the dimension of $K$ is less than $d-1$, then $n_K[x,y] = 0$ for all $[x,y] \in \mathcal{E}_X$ with probability 1 and thus $c_0, \ldots, c_{d-2}$ are vanishing. Hence

$$\mu(K) = c_{d-1}(\lambda) V_{d-1}(K) + c_d(\lambda) V_d(K).$$

If $K$ is of dimension $d-1$, then $\mu(K) = 2\nu_1(K) = c_{d-1}(\lambda) V_{d-1}(K)$. Suppose that $P$ is a polytope with facets $F \in \mathcal{F}(P)$ and nonempty interior. Then the valuation $\mu(\cdot)$ can be written as

$$\mu(P) = \sum_{F \in \mathcal{F}(P)} \mathbb{E} \sum_{[x,y] \in \mathcal{E}, [x,y] \cap F \neq \varnothing} (f(X,x) + f(X,y))$$
$$= \sum_{F \in \mathcal{F}(P)} \nu_1(F)$$
$$= c_{d-1}(\lambda) V_{d-1}(P)$$

and thus $c_d = 0$ which proves (2).

For the second claim of the theorem the scaling property of the Poisson process is used. Write $\mu_\lambda(K) = \mathbb{E} \sum_{[x,y] \in \mathcal{E}_X} (f(X,x) + f(X,y)) n_K[x,y]$ to emphasize the dependence on the intensity $\lambda$ of the point process $X$. For $t > 0$ replace $x, y$ by $x/t, y/t$. Then

$$\mu_\lambda(K) = \mathbb{E} \sum_{[x,y] \in \mathcal{E}_{tX}} t^{-\alpha}(f(tX,x) + f(tX,y)) n_{tK}[x,y] = t^{-\alpha} \mu_{t^{-d}\lambda}(tK),$$

and, together with (2), we obtain

$$c_f(\lambda) = t^{d-1-\alpha} c_f(t^{-d}\lambda) = \lambda^{(d-\alpha-1)/d} c_f(1). \qquad \square$$

**3. A jackknife estimate for the variance of functionals of a Poisson point process.** To get an estimate for the variance of a function $S(X)$ we rewrite the Efron–Stein jackknife inequality [2] (see also Efron [1] and Hall [6]). In the form we are interested in this is possible if there are no far-reaching dependencies. This is made precise in the following assumptions:

For a locally finite subset $Y \in \mathbb{R}^d$ we call $R(Y)$ the *radius of influence* of a function $S(Y)$, if there is a function $f:\mathbb{R} \to \mathbb{R}$ such that for arbitrary locally finite point sets $D \subset \mathbb{R}^d$ and $D^\circ \subset (\mathbb{R}^d \setminus R(Y)B^d)$, we have:

(A1) $S(Y) = S(Y \cup D^\circ)$;
(A2) $|S(Y) - S(Y \cup D)| \leq f(R(Y))$.

Hence the influence on $S(\cdot)$ of additional points can be estimated by the function $f(R(Y))$, and additional points outside $R(Y)B^d$ are negligible. This notion is close to the notion of *stabilization* used in previous work, for example, in Penrose [14], where stabilization at the origin refers to condition (A1) whereas here we need in addition bounds on the costs of adding points close to the origin.

We are interested in the case when $Y$ is the realization of a Poisson point process $X$ and the moments of $f(R(X))$ are bounded.

THEOREM 4. *Let $X$ (resp. $X^+$), be a Poisson point process of intensity $\lambda$, resp. $\lambda(1 + \frac{1}{m})$. Let $S:\mathbf{N} \to \mathbb{R}$ be a measurable function on the space of locally finite point sets in $\mathbb{R}^d$, and let $R(X)$ be a radius of influence of $S(X)$. If $\mathbb{E}(f(R(X))^2 R(X)^{2d})$ exists, then*

$$\operatorname{Var} S(X) \leq \lim_{m \to \infty} \mathbb{E} \sum_{x \in X^+} (S(X^+ \setminus \{x\}) - S(X^+))^2.$$

PROOF. We start with recalling the Efron–Stein jackknife inequality in its usual form. Let $Y_i$ be independent identically distributed random variables defined on some probability space, $i = 1, \ldots, m+1$. We write $Y^{(i)}$ for $(Y_1, \ldots, Y_{i-1}, Y_{i+1}, \ldots, Y_{m+1})$. If $S(Y_1, \ldots, Y_m)$ is any real symmetric function of $m$ random variables, an estimate for the expectation of $S(\cdot)$ is given by

$$\bar{S} = \frac{1}{m+1} \sum_{i=1}^{m+1} S(Y^{(i)}).$$

The Efron–Stein jackknife inequality then says that the natural estimate for the variance $\sum (S(Y^{(i)}) - \bar{S})^2$ overestimates the real variance:

$$\operatorname{Var} S \leq \mathbb{E} \sum_{i=1}^{m+1} (S(Y^{(i)}) - \bar{S})^2.$$

Since the right-hand side increases if we replace the mean $\bar S$ by any other function $T = T(Y_1, \ldots, Y_{m+1})$, we also have

$$\text{Var } S \leq \mathbb{E} \sum_{i=1}^{m+1} (S(Y^{(i)}) - T)^2. \tag{3}$$

Let $X_1, \ldots, X_{m+1}$ be independent Poisson point processes in $\mathbb{R}^d$ of intensity $\lambda/m$, set $X^{(i)} = \bigcup_{j=1,\ldots,i-1,i+1,\ldots,m+1} X_j$, $X = X^{(m+1)}$, and $X^+ = \bigcup_{j=1}^{m+1} X_j$ which are Poisson point processes in $\mathbb{R}^d$ of intensity $\lambda$, or $\lambda(1 + \frac{1}{m})$ respectively. Since by assumption $S(X_1, \ldots, X_m) = S(X_1 \cup \cdots \cup X_m)$ is a symmetric function in the $X_i$, the Efron–Stein jackknife inequality (3) with $T = S(X^+)$ tells us that

$$\text{Var } S(X) \leq \mathbb{E} \sum_{i=1}^{m+1} (S(X^{(i)}) - S(X^+))^2$$

$$= \mathbb{E} \sum_{i=1}^{m+1} \mathbb{E}^{(i)}(S(X^{(i)}) - S(X^+))^2,$$

where $\mathbb{E}^{(i)}(\cdot)$ abbreviates $\mathbb{E}(\cdot | X^{(i)})$.

For the next step fix $i$ and denote the radius of influence of $S(X^{(i)})$ by $R(X^{(i)}) = R$. Observe that this implies that $R$ is independent of $X_i$. So we may apply conditions (A1) and (A2) with $Y = X^{(i)}$, $D^\circ = X_i \setminus RB^d$, $D = X_i$.

If $m$ is large, then with high probability at most one point $x \in X_i$ is in $RB^d$ and thus may have influence on $S(X^+)$. This is made precise in the following. We decompose the expectation according to the value of $X_i(RB^d)$:

$$\mathbb{E}^{(i)}(S(X^{(i)}) - S(X^+))^2$$

$$= \sum_{n=0}^{\infty} \mathbb{E}^{(i)}((S(X^{(i)}) - S(X^+))^2 I(X_i(RB^d) = n)).$$

For $X_i(RB^d) \in \{0, 1\}$ we use (A1) with $D^\circ = X_i \setminus RB^d$ and obtain

$$\sum_{n=0}^{1} \mathbb{E}^{(i)} \left( \sum_{x \in X_i \cap RB^d} (S(X^+ \setminus \{x\}) - S(X^+))^2 I(X_i(RB^d) = n) \right)$$

$$\leq \sum_{n=0}^{\infty} \mathbb{E}^{(i)} \left( \sum_{x \in X_i} (S(X^+ \setminus \{x\}) - S(X^+))^2 I(X_i(RB^d) = n) \right)$$

$$= \mathbb{E}^{(i)} \sum_{x \in X_i} (S(X^+ \setminus \{x\}) - S(X^+))^2.$$

For $X_i(RB^d) \geq 2$ we use (A2) with $D = X_i$, the estimate $\sum_{n=2}^{\infty} \frac{\mu^n}{n!} e^{-\mu} \leq \frac{\mu^2}{2}$, and obtain

$$\sum_{n=2}^{\infty} \mathbb{E}^{(i)}((S(X^{(i)}) - S(X^+))^2 I(X_i(RB^d) = n)) \leq \sum_{n=2}^{\infty} f(R)^2 \mathbb{P}(X_i(RB^d) = n)$$

$$\leq f(R)^2 \frac{\lambda^2 V(RB^d)^2}{2m^2}$$

since the intensity of $X_i$ equals $\frac{\lambda}{m}$.

Combining our results, summing over $i = 1, \ldots, m+1$, and using that the radii of influence $R(X^{(i)})$ are identically distributed, gives

$$(4) \quad \operatorname{Var} S(X) \leq \mathbb{E} \sum_{x \in X^+} (S(X^+ \setminus \{x\}) - S(X^+))^2 + \frac{\lambda^2}{m} \mathbb{E}(f(R)^2 V(RB^d)^2)$$

and thus proves the theorem. □

In the next section we use Theorem 4 for functionals with moments continuous in the intensity of $X$. From (4) we obtain in this case the following corollary.

COROLLARY 5. *Let $X$ be a Poisson point process. Let $S : \mathbf{N} \to \mathbb{R}$ be a measurable function on the space of locally finite point sets in $\mathbb{R}^d$, and let $R(X)$ be a radius of influence of $S(X)$. If $\mathbb{E} S(X)$, $\mathbb{E} S^2(X)$ are continuous in $\lambda$, and if $\mathbb{E}(f(R(X))^2 R(X)^{2d})$ exists, then*

$$\operatorname{Var} S(X) \leq \mathbb{E} \sum_{x \in X} (S(X \setminus \{x\}) - S(X))^2.$$

We want to remark that the Slivnyak–Mecke formula for a Poisson point process allows to rewrite our theorem in the following way:

$$\operatorname{Var} S(X) \leq \lambda \int_{\mathbb{R}^d} \lim_{m \to \infty} \mathbb{E}(S(X^+) - S(X^+ \cup \{x\}))^2 \, dx.$$

We conjecture that the following more general theorem holds:

CONJECTURE. *Let $X$ be a Poisson point process. For any measurable function $S : \mathbf{N} \to \mathbb{R}$ on the space of locally finite point sets in $\mathbb{R}^d$ we have*

$$\operatorname{Var} S(X) \leq \mathbb{E} \sum_{x \in X} (S(X \setminus \{x\}) - S(X))^2.$$

**4. Volume difference.** In this section we are interested in the difference of the volume of $v_X(K)$ and $K$. We state the mean value and prove an estimate for the variance. The large deviation inequality is proved in Section 6.

It can easily be shown that

$$\mathbb{E}V(v_X(K)) = V(K). \tag{5}$$

This follows either from Campbell's theorem (see, e.g., the book by Schneider and Weil [17], page 128) or using Hadwiger's characterization theorem and an argument similar to that of Theorem 3. Formula (5) holds for all Borel sets without any convexity assumptions.

To get a bound on the variance of $V(v_X(K))$ we use the Efron–Stein jackknife inequality in Corollary 5. This states that

$$\operatorname{Var} V(v_X(K)) \leq \mathbb{E} \sum_{x \in X} (V(v_{X \setminus \{x\}}(K)) - V(v_X(K)))^2$$

if for some radius $R(X)$ of influence the moment $\mathbb{E}(f(R(X))^2 R(X)^{2d})$ exists. (Observe that the moments of the functional we are interested in are continuous in $\lambda$.)

Thus we have to estimate the volume of those Voronoi cells with centers $x \in X$, which partly may change from exterior points to interior points or vice versa if $x$ is removed. Assume that $X \cap \partial K$ is empty which happens with probability 1.

If for $x \in X \cap K$ all neighbors of the Voronoi cell $v_X(x)$ are also contained in $K$, that is, if for all $[x,y] \in \mathcal{E}_X$ we have $y \in K$, then $v_{X \setminus \{x\}}(K) = v_X(K)$. The same argument applies if the point $x$ and all its neighbors are outside $K$. Hence of interest are those points $x \in X$ such that there exists an edge $[x,y] \in \mathcal{E}_X$ with $[x,y] \cap \partial K \neq \varnothing$ in which case

$$|V(v_{X \setminus \{x\}}(K)) - V(v_X(K))| \leq V(v_X(x)).$$

Defining $n_K[x,y]$ as in Section 2 and noting that $n_K[x,y] \geq \mathbf{1}([x,y] \cap \partial K \neq \varnothing)$ we thus see that

$$\operatorname{Var} V(v_X(K)) \leq \mathbb{E} \sum_{[x,y] \in \mathcal{E}_X} n_K[x,y](V(v_X(x))^2 + V(v_X(y))^2).$$

By Theorem 3 with $\alpha = 2d$ we immediately obtain

$$\operatorname{Var} V(v_X(K)) \leq c_d \lambda^{-1-1/d} S(K)$$

which is the variance estimate of Theorem 1.

It remains to define the radius of influence $R(X)$, and to show that $\mathbb{E}(f(R(X))^2 R(X)^{2d})$ exists. Define a (random) number $R' = R'(X)$ as the smallest number fulfilling

$$\bigcup_{z \in K} v_X(z) \subset R' B^d \tag{6}$$

and let $R(X) = 3R'$. [Recall that $v_X(K)$ is the union of all Voronoi cells $v_X(x)$ with nucleus $x \in X \cap K$, whereas here we estimate the influence of all $z \in K$.]

As for assumption (A1) we have to show that any point set $D°$ which does not meet $3R'B^d$ has no influence. Indeed, if $v_X(K) \neq v_{X \cup D°}(K)$, then there are points $x \in X \cap K$ and $y \in D°$ with $[x,y] \in \mathcal{E}_{X \cup D°}$. Hence the Voronoi cells of $x$ and $y$ would have points in common. This is impossible since by the definition of $R'$ we have $v_X(x) \subset R'B^d$, but the midplane between $x$ and $y$ does not meet the ball $R'B^d$.

As for assumption (A2), it follows from (6) that for any point set $D$

$$v_{X \cup D}(K) \subset v_{X \cup (D \cap K)}(K) \subset R'B^d.$$

Thus the difference between $v_X(K)$ and $v_{X \cup D}(K)$ is bounded by the volume of $R'B^d$ and assumption (A2) is fulfilled with $f(R(X)) = V(R'B^d)$.

Finally we have to show that $\mathbb{E}(R'^{4d})$ is finite. Denote by $R_K$ the smallest radius such that $K \subset R_K B^d$. By definition, if $R' \geq r$ for some $r \geq R_K + \sqrt{d}$, then there is a point $y \in rS^{d-1}$ with $d(y, \partial K) \leq d(y, X)$, that is, $X(B(y, (r - R_K))) = 0$. We cover the ball $rB^d$ by $2^d r^d$ disjoint cubes $C_i$ of sidelength 1 with center $z_i$ and obtain

$$\mathbb{P}(R' \geq r) \leq \sum_{i=1}^{2^d r^d} \mathbb{P}(\exists y \in rS^{d-1} \cap C_i : X(B(y, (r - R_K))) = 0)$$

$$\leq \sum_{i=1}^{2^d r^d} \mathbb{P}(\exists y \in C_i : X(B(z_i, (r - \sqrt{d} - R_K))) = 0)$$

$$\leq 2^d r^d e^{-\lambda \kappa_d (r - \sqrt{d} - R_K)^d}.$$

Thus all moments of $R'$ exist and are finite.

**5. Symmetric difference metric.** In this section we investigate the volume of the symmetric difference of $v_X(K)$ and $K$,

$$V(K \triangle v_X(K)) = V(K \setminus v_X(K)) + V(v_X(K) \setminus K).$$

We determine the expectation (Theorem 6), and prove an estimate for the variance. The large deviation inequality is proved in Section 6.

First we show that

$$V(K \triangle v_X(K)) = c'_{d-1} S(K) \lambda^{-1/d} + o(\lambda^{-1/d}).$$

We start with the volume of $v_X(K) \setminus K$. The Slivnyak–Mecke formula gives for $x \in \mathbb{R}^d \setminus K$

$$\mathbb{P}(x \in v_X(K)) = \mathbb{P}(\exists y \in X \cap K : x \in v_X(y))$$

$$= \mathbb{E} \sum_{y \in X \cap K} \mathbf{1}(x \in v_X(y))$$

(7)
$$= \lambda \int_K \mathbb{P}(x \in v_{X \cup \{y\}}(y)) \, dy$$

$$= \lambda \int_K e^{-\lambda V(B(x,d(x,y)))} \, dy$$

since $x \in v_{X \cup \{y\}}(y)$ if the intersection of $X$ with the ball of radius $d(x,y)$ centered at $x$ is empty. Precisely the same argument shows that for $x \in K$

(8) $\mathbb{P}(x \notin v_X(K)) = \mathbb{E} \sum_{y \in X \setminus K} \mathbf{1}(x \in v_X(y)) = \lambda \int_{\mathbb{R}^d \setminus K} e^{-\lambda V(B(x,d(x,y)))} \, dy.$

Combining (7) and (8) we obtain

$$\mathbb{E} V(K \triangle v_X(K)) = \mathbb{E} \int_{\mathbb{R}^d} \mathbf{1}(x \in K \triangle v_X(K)) \, dx$$

$$= \int_{\mathbb{R}^d \setminus K} \mathbb{P}(x \in v_X(K)) \, dx + \int_K \mathbb{P}(x \notin v_X(K)) \, dx$$

$$= 2\lambda \int_{\mathbb{R}^d \setminus K} \int_K e^{-\lambda \kappa_d \|y - x\|^d} \, dy \, dx.$$

We use the Blaschke–Petkantschin formula (see, e.g., [16]) which transforms the integration of the tuple $(x,y)$ with respect to Lebesgue measure into integration of $(x,y)$ with respect to the (one-dimensional) Lebesgue measure on the line $E$ which is the affine hull of the two points, and then integrate with respect to the set $\mathcal{E}_1^d$ of all lines in $\mathbb{R}^d$ using the normalized Haar measure $\nu$ on the set of all lines:

$$\mathbb{E} V(K \triangle v_X(K)) = d\lambda \kappa_d \int_{\mathcal{E}_1^d} \int_{E \setminus K} \int_{E \cap K} e^{-\lambda \kappa_d \|y-x\|^d} \|y-x\|^{d-1} \, dy \, dx \, d\nu(E).$$

Identify $E$ with $\mathbb{R}$ and $E \cap K$ with the interval $[0,l]$ of length $l = l(E)$. If $l > 0$, we obtain for the inner integrations

$$\int_{\mathbb{R} \setminus [0,l]} \int_0^l e^{-\lambda \kappa_d |x-y|^d} |x-y|^{d-1} \, dy \, dx = 2 \int_0^l \int_l^\infty e^{-\lambda \kappa_d (x-y)^d} (x-y)^{d-1} \, dx \, dy$$

$$= \frac{2}{d}(\lambda \kappa_d)^{-1} \int_0^l e^{-\lambda \kappa_d y^d} \, dy$$

$$= \frac{2}{d^2}(\lambda \kappa_d)^{-1-1/d} \int_0^{\lambda \kappa_d l^d} e^{-s} s^{1/d-1} \, ds$$

$$= \frac{2}{d^2}(\lambda \kappa_d)^{-1-1/d} \Gamma\left(\frac{1}{d}\right)(1 - \delta(\lambda, E)),$$

where
$$0 \leq \Gamma\left(\frac{1}{d}\right)\delta(\lambda, E) = \int_{\lambda\kappa_d l^d}^{\infty} e^{-s} s^{(1/d)-1} \, ds \leq \Gamma\left(\frac{1}{d}\right) e^{-\lambda\kappa_d l^d}$$

since $\frac{1}{d} - 1 < 0$. Thus we have

$$\mathbb{E}V(K \triangle v_X(K)) = \frac{2}{d}(\lambda\kappa_d)^{-1/d}\Gamma\left(\frac{1}{d}\right) \int_{\mathcal{E}_1^d} \mathbf{1}(E \cap K \neq \varnothing)(1 - \delta(\lambda, E)) \, d\nu(E).$$

For the main term we obtain by Cauchy's surface area formula

$$\frac{2}{d}(\lambda\kappa_d)^{-1/d}\Gamma\left(\frac{1}{d}\right) \int_{\mathcal{E}_1^d} \mathbf{1}(E \cap K \neq \varnothing) \, d\nu(E)$$

$$= \frac{2}{d}(\lambda\kappa_d)^{-1/d}\kappa_{d-1}\Gamma\left(\frac{1}{d}\right) S(K).$$

To estimate the error term $\int \delta(\lambda, E) \, d\nu(E)$ assume that the origin of the coordinate system is chosen in such a way that

$$r(K)B^d \subset K,$$

where $r(K)$ is the inradius of $K$. Parametrize the line $E$ by $E = tu + y$, $t \in \mathbb{R}$, where $u \in S^{d-1}$ is the direction of $E$ and $y \in u^\perp$. The measure $\nu$ decomposes into the uniform distribution $\omega$ on the sphere $S^{d-1}$, and for $u \in S^{d-1}$, into Lebesgue measure in the hyperplane $u^\perp$. If the line $E$ meets $K$, then the point $y$ is contained in the projection $K|_{u^\perp}$ of the set $K$ onto $u^\perp$:

$$\int_{\mathcal{E}_1^d} \mathbf{1}(E \cap K \neq \varnothing)\delta(\lambda, E) \, d\nu(E) = \int_{S^{d-1}} \int_{K|_{u^\perp}} \delta(\lambda, E) \, dy \, d\omega(u).$$

We introduce polar coordinates $y = rv$, where $v$ is integrated with respect to Lebesgue measure $\sigma$ on $S^{d-1} \cap u^\perp$. Denote by $\rho(v) = \rho_{K|_{u^\perp}}(v)$ the radial function of $K|_{u^\perp}$ in direction $v$. Because $K$ is a convex set, we have that for fixed $u$ and $v$ the chord length $l(rv, u)$ is a concave function in $r$ which vanishes at the boundary of $K|_{u^\perp}$. Hence

$$l(rv, u) \geq l(0, u)\left(1 - \frac{r}{\rho(v)}\right) \geq 2r(K)\left(1 - \frac{r}{\rho(v)}\right)$$

for $0 \leq r \leq \rho(v)$. This yields

$$\int_{K|_{u^\perp}} e^{-\lambda\kappa_d l^d} \, dy \leq \int_{S^{d-1} \cap u^\perp} \int_0^{\rho(v)} e^{-\lambda\kappa_d 2^d r(K)^d (1 - r/\rho(v))^d} r^{d-2} \, dr \, d\sigma(v)$$

$$\leq \frac{1}{2d}(\lambda\kappa_d)^{-1/d} r(K)^{-1} \int_{S^{d-1} \cap u^\perp} \rho(v)^{d-1} \, d\sigma(v) \int_0^{\infty} e^{-s} s^{1/d-1} \, ds$$

$$\leq \frac{d-1}{2d}(\lambda\kappa_d)^{-1/d}\Gamma\left(\frac{1}{d}\right) r(K)^{-1} V_{d-1}(K|_{u^\perp}).$$

Using Cauchy's surface area formula again gives

$$\int_{\mathcal{E}_1^D} \mathbf{1}(E \cap K \neq \varnothing)\delta(\lambda, E)\, d\nu(E) \leq \frac{d-1}{2d}(\lambda\kappa_d)^{-1/d}\Gamma\left(\frac{1}{d}\right)r(K)^{-1}\kappa_{d-1}S(K).$$

We summarize our results:

THEOREM 6. *If $K \in \mathcal{K}^d$, then*

$$\mathbb{E}V(K \triangle v_X(K)) = \frac{2}{d}(\lambda\kappa_d)^{-1/d}\kappa_{d-1}\Gamma\left(\frac{1}{d}\right)S(K)(1 - \lambda^{-1/d}\Delta),$$

*where* $0 \leq \Delta \leq \frac{d-1}{2d}\kappa_d^{-1/d}\Gamma(\frac{1}{d})r(K)^{-1}$.

The same arguments which led to the bound on the variance of $V(v_X(K))$ will yield a bound on the variance of $V(K \triangle v_X(K))$. We use again the Efron–Stein jackknife inequality proved in Theorem 4 showing that

$$\operatorname{Var} V(K \triangle v_X(K)) \leq \mathbb{E}\sum_{x \in X}(V(v_{X\setminus\{x\}}(K)\triangle K) - V(v_X(K)\triangle K))^2$$

where the radius of influence $R(X)$ is defined precisely as in (6). Hence we already know that $\mathbb{E}(f(R(X))^2 R(X)^{2d})$ exists.

As in Section 4, of interest are those points $x \in X$ for which there exists an edge $[x, y] \in \mathcal{E}_X$ with $[x, y] \cap \partial K \neq \varnothing$. In this case

$$|V(v_{X\setminus\{x\}}(K)\triangle K) - V(v_X(K)\triangle K)| \leq V(v_X(x)).$$

Thus we obtain

$$\operatorname{Var} V(K \triangle v_X(K)) \leq \mathbb{E}\sum_{[x,y] \in \mathcal{E}_X} n_K[x,y](V(v_X(x))^2 + V(v_X(y))^2)$$

$$\leq c_f \lambda^{-1-1/d}S(K)$$

which is the variance estimate of Theorem 2.

**6. Large deviation inequalities.** In this section we prove the large deviation inequalities of Theorems 1 and 2. The essential tool is Azuma's inequality, in particular the method of *uniformly difference-bounded functions* used by McDiarmid [10].

A function $f: \Omega_1 \times \cdots \times \Omega_m \to \mathbb{R}$ is called *uniformly difference-bounded by $b$* if the following holds: for any $(y_1, \ldots, y_m) \in \Omega_1 \times \cdots \times \Omega_m$, and for any $k$ and any $y_k' \in \Omega_k$ we have

$$|f(y_1, \ldots, y_k, \ldots, y_m) - f(y_1, \ldots, y_k', \ldots, y_m)| \leq b.$$

Let $Y_1, \ldots, Y_m$ be independent random variables with $Y_k \in \Omega_k$. Set $\zeta = f(Y_1, \ldots, Y_m)$ with $f$ uniformly difference-bounded by $b$. Then McDiarmid's bounded difference inequality says that, for any $t$,

$$\mathbb{P}(|\zeta - \mathbb{E}\zeta| \geq t) \leq 2e^{-2t^2/(mb^2)}.$$

To define the random variables $Y_1, \ldots, Y_m$ we need some preparations. Dissect $\mathbb{R}^d$ into cubes $C_i$ of diameter $\delta$ having pairwise disjoint interior. Define $\delta$ such that

$$V(C_i) = d^{-d/2}\delta^d = k\lambda^{-1}\ln\lambda$$

with $k \geq 2$. Assume that the cubes are numbered in such a way that for $i = 1, \ldots, m_\delta$ the cubes $C_i$ have nonempty intersection with $\partial K + 3\delta B^d$, and that $C_i$ is disjoint from $\partial K + 3\delta B^d$ for $i > m_\delta$. Since for $i = 1, \ldots, m_\delta$ the cubes $C_i$ are contained in $\partial K + 4\delta B^d$, we see that

$$\sum_{i=1}^{m_\delta} V(C_i) = d^{-d/2}\delta^d m_\delta \leq V(\partial K + 4\delta B^d).$$

Note that there is a $\Theta(K) \in (0, 1]$ such that for $\varepsilon \leq \Theta(K)$ we have

(9) $$V(\partial K + \varepsilon B^d) \leq 4\varepsilon S(K).$$

Thus for $4\delta \leq \Theta(K)$ we obtain

(10) $$m_\delta \leq 16 d^{d/2} S(K) \delta^{-d+1}.$$

For $\lambda$ large, in each of the cubes $C_i$ at least one point of the Poisson point process is contained with high probability. To make this precise denote by $A$ the event that for all $i = 1, \ldots, m_\delta$ we have $X(C_i) \geq 1$. Since $X \cap C_i$ is Poisson distributed we have

(11) $$\mathbb{P}(A^c) = 1 - (1 - e^{-\lambda V(C_i)})^{m_\delta} \leq m_\delta e^{-\lambda V(C_i)} = m_\delta \lambda^{-k}.$$

We assume in the following that each cube $C_i$, $i = 1, \ldots, m_\delta$, contains at least one point. This implies that if $x \in X$ has distance at most $\delta$ to $\partial K$, then

(12) $$v_X(x) \subset B(x, \delta).$$

To prove this inclusion observe that for any point $y \in C_i$, $i = 1, \ldots, m_\delta$, the distance of $y$ to one of the points $X \cap C_i$ is at most $\delta$, the diameter of the cube $C_i$. Thus if $y$ is contained in some Voronoi cell $v_X(\tilde{x})$ with $\tilde{x} \in X$, then $\tilde{x}$ is the nearest point of $X$ to $y$, and we have

(13) $$\|y - \tilde{x}\| \leq \delta.$$

Now let $x \in X$ have distance at most $\delta$ to the boundary of $K$. Assume

$$y \in B(x, 2\delta) \setminus B(x, \delta).$$

Then the distance of $y$ to the boundary of $K$ is at most $3\delta$ and thus $y$ is contained in some $C_i$, $i = 1, \ldots, m_\delta$. By (13) the distance from $y$ to the nearest point of $X \cap C_i$ is at most $\delta$ and hence $\|y - x\| > \delta$ implies that $y \notin v_X(x)$. This proves (12) since the Voronoi cell $v_X(x)$ does not meet $B(x, 2\delta) \setminus B(x, \delta)$ and is connected.

Since all Voronoi cells meeting the boundary of $K$ have circumradius at most $\delta$, they are contained in $\partial K + 2\delta B^d$. By (13) the centers of neighbors cells, having boundary points $y$ in common with these cells, have distance at most $3\delta$ to $\partial K$ and thus are contained in $C_i$, $i = 1, \ldots, m_\delta$. In other words, the set of all Voronoi cells meeting the boundary of $K$ only depends on $X \cap C_i$, $i = 1, \ldots, m_\delta$, and is independent of all points of $X$ outside the cubes $C_i$, $i = 1, \ldots, m_\delta$.

For $i = 1, \ldots, m_\delta$ define the random points $Y_i$ by $Y_i = X \cap C_i$. If $\zeta = f(X)$ is a function depending only on those Voronoi cells meeting the boundary of $K$, then $\zeta$ depends only on $Y_i$, $\zeta = f(X) = f(Y_1, \ldots, Y_{m_\delta})$. In the cases we are interested in, we have either

$$\zeta = f(Y_1, \ldots, Y_{m_\delta}) = V(v_X(K)) - V(K)$$
$$= \sum_{i=1}^{m_\delta} \left( \sum_{x \in Y_i \cap K} V(v_X(x) \setminus K) - \sum_{x \in Y_i \setminus K} V(v_X(x) \cap K) \right)$$

or

$$\zeta = f(Y_1, \ldots, Y_{m_\delta}) = V(K \triangle v_X(K))$$
$$= \sum_{i=1}^{m_\delta} \left( \sum_{x \in Y_i \cap K} V(v_X(x) \setminus K) + \sum_{x \in Y_i \setminus K} V(v_X(x) \cap K) \right).$$

In both cases it follows from (12) that, replacing $Y_i$ by some nonempty finite subset $Y_i' \subset C_i$, we have

$$|f(\ldots, Y_i, \ldots) - f(\ldots, Y_i', \ldots)| \leq V(C_i + \delta B^d)$$
$$\leq 3^d \delta^d$$

and thus $b = 3^d \delta^d$ is the required difference bound. Now McDiarmid's theorem tells us that

$$\mathbb{P}(|\zeta - \mathbb{E}(\zeta|A)| \geq t|A) \leq 2e^{-2t^2/(m_\delta b^2)}.$$

It follows from $\mathbb{P}(\cdot) \leq \mathbb{P}(\cdot|A) + \mathbb{P}(A^c)$ and from (11) that

$$\mathbb{P}(|\zeta - \mathbb{E}(\zeta|A)| \geq t) \leq 2e^{-2t^2/(m_\delta b^2)} + m_\delta \lambda^{-k}.$$

In the last step we replace $\mathbb{E}(\zeta|A)$ by $\mathbb{E}\zeta$. We use the elementary inequality

$$|\mathbb{E}\zeta - \mathbb{E}(\zeta|A)| \leq |\mathbb{E}(\zeta \mathbf{1}(A)) - \mathbb{E}(\zeta|A)| + \mathbb{E}(\zeta \mathbf{1}(A^c))$$

$$\leq \mathbb{E}(\zeta|A)\mathbb{P}(A^c) + \sqrt{\mathbb{E}(\zeta^2)\mathbb{P}(A^c)}$$
$$\leq (\mathbb{E}(\zeta|A) + \sqrt{\mathbb{E}(\zeta^2)})\sqrt{\mathbb{P}(A^c)},$$

where the second line follows from Hölder's inequality. Since, conditioning on $A$, all Voronoi cells meeting the boundary of $K$ are contained in $\partial K + 2\delta B^d$, we have by (9)

$$\mathbb{E}(\zeta|A) \leq 8\delta S(K).$$

And for $\mathbb{E}\zeta^2$ the bounds on the expectation and variance yield immediately that

$$\mathbb{E}\zeta^2 \leq c_2(K)\lambda^{-2/d}$$

for $\lambda \geq 1$ (which follows from $\delta \leq 1$). Thus

$$2|\mathbb{E}\zeta - \mathbb{E}(\zeta|A)|^2 \; m_\delta^{-1}b^{-2} \leq c_3(K)k^2\lambda^{-k-2/d}\delta^{-2d}$$
$$\leq c_4(K)\lambda^{-2/d}$$

since $\delta \leq 1$ and $k \geq 2$.

Define $x_+ = \max(0,x)$. Using the inequality $2(t-s)_+^2 \geq t^2 - 2s^2$ and (10) we obtain

$$\mathbb{P}(|\zeta - \mathbb{E}\zeta| \geq t) \leq 2e^{-2(t-|\mathbb{E}\zeta-\mathbb{E}(\zeta|A)|)_+^2/(m_\delta b^2)} + m_\delta\lambda^{-k}$$
$$\leq 2e^{c_4(K)\lambda^{-2/d}}e^{-t^2/(m_\delta b^2)} + m_\delta\lambda^{-k}$$
$$\leq c_5(K)e^{-c_d t^2(k\ln\lambda)^{-1-1/d}\lambda^{1+1/d}S(K)^{-1}} + 16\sqrt{d}S(K)\lambda^{-k+1}$$

with $c_d = 2^{-4}3^{-2d}d^{-d-1/2}$ for $4\sqrt{d}(k\lambda^{-1}\ln\lambda)^{1/d} \leq \Theta(K)$ and any $k \geq 2$.

## REFERENCES

[1] EFRON, B. (1982). *The Jackknife, the Bootstrap and Other Resampling Plans. CBMS–NSF Regional Conference Series in Applied Mathematics* **38**. SIAM, Philadelphia, PA. MR659849
[2] EFRON, B. and STEIN, C. (1981). The jackknife estimate of variance. *Ann. Statist.* **9** 586–596. MR615434
[3] EINMAHL, J. H. J. and KHMALADZE, E. V. (2001). The two-sample problem in $\mathbb{R}^m$ and measure-valued martingales. In *State of the Art in Probability and Statistics (Leiden, 1999). IMS Lecture Notes—Monograph Series* **36** 434–463. IMS, Beachwood, OH. MR1836574
[4] GRAF, S. and LUSCHGY, H. (2000). *Foundations of Quantization for Probability Distributions. Lecture Notes in Mathematics* **1730**. Springer, Berlin. MR1764176
[5] HADWIGER, H. (1957). *Vorlesungen über Inhalt, Oberfläche und Isoperimetrie.* Springer, Berlin. MR0102775
[6] HALL, P. (1992). *The Bootstrap and Edgeworth Expansion.* Springer, New York. MR1145237

[7] KHMALADZE, E. and TORONJADZE, N. (2001). On the almost sure coverage property of Voronoi tessellation: The $\mathbb{R}^1$ case. *Adv. in Appl. Probab.* **33** 756–764. MR1875777

[8] KLAIN, D. A. (1995). A short proof of Hadwiger's characterization theorem. *Mathematika* **42** 329–339. MR1376731

[9] KLAIN, D. A. and ROTA, G.-C. (1997). *Introduction to Geometric Probability*. Cambridge Univ. Press. MR1608265

[10] MCDIARMID, C. (1989). On the method of bounded differences. In *Surveys in Combinatorics, 1989 (Norwich, 1989). London Mathematical Society Lecture Note Series* **141** 148–188. Cambridge Univ. Press. MR1036755

[11] MCMULLEN, P. (1993). Valuations and dissections. In *Handbook of Convex Geometry, Vol. B* (P. M. GRUBER AND J. MILLS, EDS.) 933–988. North-Holland, Amsterdam. MR1243000

[12] MCMULLEN, P. and SCHNEIDER, R. (1983). Valuations on convex bodies. In *Convexity and Its Applications* (P. M. GRUBER AND J. MILLS, EDS.) 170–247. Birkhäuser, Basel. MR731112

[13] MØLLER, J. (1994). *Lectures on Random Voronoĭ Tessellations. Lecture Notes in Statistics* **87**. Springer, New York. MR1295245

[14] PENROSE, M. D. (2007). Laws of large numbers in stochastic geometry with statistical applications. *Bernoulli* **13** 1124–1150. MR2364229

[15] SCHNEIDER, R. (1993). *Convex Bodies: The Brunn–Minkowski Theory. Encyclopedia of Mathematics and Its Applications* **44**. Cambridge Univ. Press. MR1216521

[16] SCHNEIDER, R. and WEIL, W. (1992). *Integralgeometrie*. Teubner, Stuttgart. MR1203777

[17] SCHNEIDER, R. and WEIL, W. (2000). *Stochastische Geometrie*. Teubner, Stuttgart. MR1794753

[18] STOYAN, D., KENDALL, W. S. and MECKE, J. (1987). *Stochastic Geometry and Its Applications*. Wiley, Chichester. MR895588

INSTITUTE OF DISCRETE MATHEMATICS
AND GEOMETRY
VIENNA UNIVERSITY OF TECHNOLOGY
WIEDNER HAUPTSTRASSE 8
1040 VIENNA
AUSTRIA
E-MAIL: matthias.reitzner@tuwien.ac.at